# Geometric Decomposition-Based Formulation for Time Derivatives of Instantaneous Impact Point

Byeong-Un Jo[1] and Jaemyung Ahn[2]

**Abstract:** A new analytic formulation to express the time derivatives of the instantaneous impact point (IIP) of a rocket is proposed. The geometric relationship on a plane tangential to the IIP is utilized to decompose the inertial IIP rate vector into the downrange and crossrange components, and a systematic procedure to determine the component values is presented. The new formulation shows significant advantages over the existing formulation such that the procedure and final expressions for the IIP derivatives are easy to understand and more compact. The validity of the proposed formulation was demonstrated through numerical simulation.

## Introduction

An instantaneous impact point (IIP) is the point of intersection between the orbit of a rocket determined by its current position and velocity and the surface of the Earth, which can be interpreted as the point of impact on the ground with free-fall flight (or immediate cut-off of propelled flight) assumption of the rocket. The IIP is a very important information for flight safety operations during the launch campaign. It should be calculated in real time (at the launch support facility on the ground or at the on-board computer of the rocket) and monitored continuously in order to detect potential causalities or economic damage due to the falling of rocket debris and take responsive actions (e.g., command flight termination) in a timely manner. A number of studies on IIP calculation/application with various computational procedures, modeling assumptions, and correction methodologies have been published including the algorithm with flat-Earth assumption subject to constant gravity (Montenbruck et al. 2002; Montenbruck and Markgraf 2004), the Keplerian IIP in a spherical/elliptical Earth (Ahn et al. 2000; Ahn and Roh 2012), and correction of atmospheric drag using the response surface model (Ahn and Seo 2013), which was applied in launch vehicle mission design (Yoon and Ahn 2015).

---
[1]  Graduate Research Assistant, Department of Aerospace Engineering, Korea Advanced Institute of Science and Technology (KAIST), 291 Daehak-Ro, Daejeon 34141, Republic of Korea.
[2]  Associate Professor, Department of Aerospace Engineering, Korea Advanced Institute of Science and Technology (KAIST), 291 Daehak-Ro, Daejeon 34141, Republic of Korea (corresponding author). E-mail: jaemyung.ahn@kaist.ac.kr

When there exist external accelerations with the exception of gravity – such as thrust and atmospheric drag – a change in IIP occurs, which involves the concept of "time derivatives of IIP" or "IIP rate." Studies on the computation and application of the IIP rate are very scarce. Ahn and Roh (2014) developed analytic formulations to express the time derivatives of Keplerian IIP, which were utilized by Nam et al. (2015) to improve the flight termination decision procedure used for the flight safety system of a rocket.

While the analytic formulations presented by Ahn and Roh (2014) are mathematically solid and implementable, there is still room for the improvement of their procedure and results. The procedure highly depends on the mechanical differentiation of the final expression for IIP, which lacks the physical interpretations of the parameters that appear in the formulation. In addition, the final results representing the IIP rate are not properly arranged (unnecessarily complex), and can be misleading in terms of independent direction vectors determining the overall IIP rate.

This paper proposes a new formulation to describe the time derivatives of the IIP of a rocket and address the aforementioned disadvantages of existing study. The inertial IIP rate vector is decomposed into the downrange and crossrange components based on the geometry on a plane tangential to IIP, and a systematic procedure to determine the component values is presented. The new formulation significantly streamlines the computational procedure with the improved clarity provided by the physical interpretation. The proposed formulation has been verified through comparison with an existing algorithm through numerical simulation.

# Review: Instantaneous Impact Point and its Time Derivatives

## *Analytic Expression for Keplerian IIP*

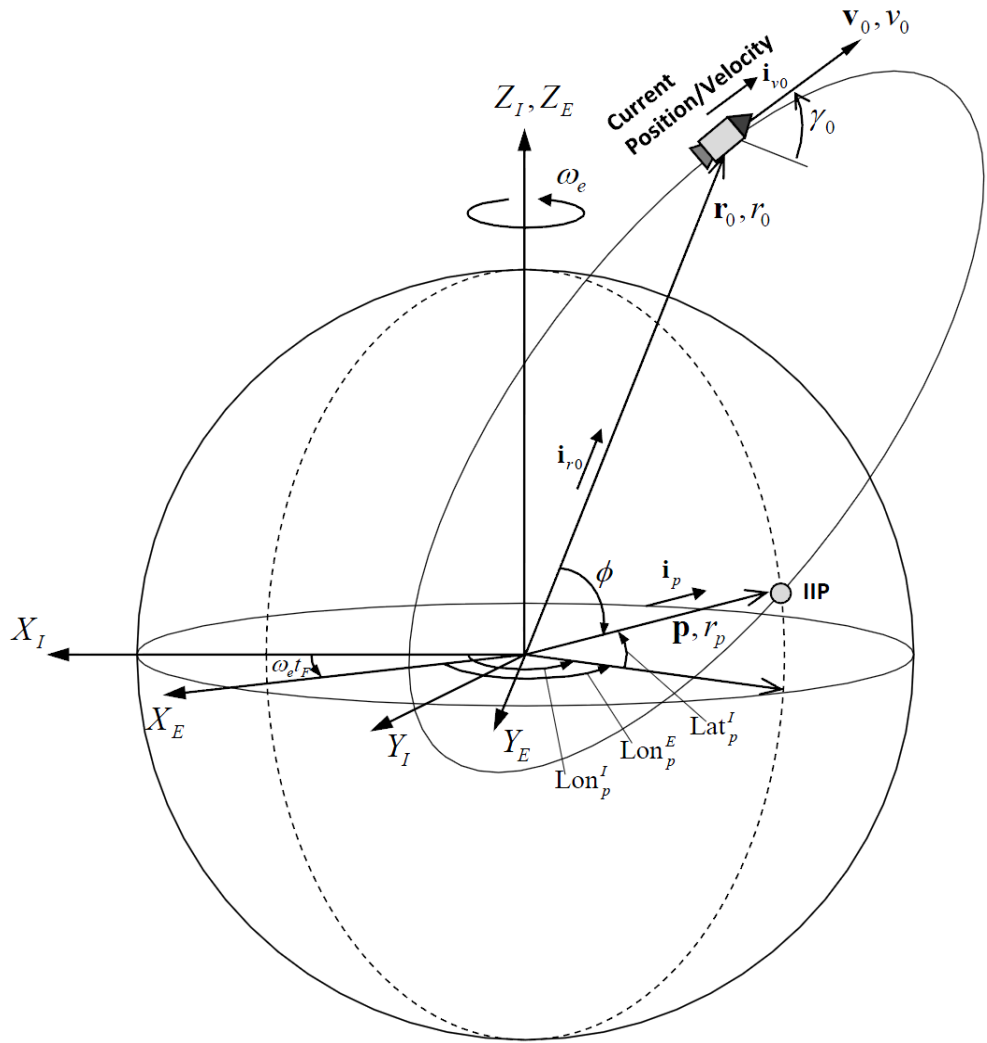

**Fig. 1. Geometry for calculating the IIP of a rocket**

Fig. 1 shows the geometry used for calculating the IIP ($\mathbf{p}$) of a rocket with given initial position ($\mathbf{r}_0$) and velocity ($\mathbf{v}_0$). The translational motion of a rocket in an Earth-centered inertial (ECI) coordinate frame can be expressed as

$$\ddot{\mathbf{r}} = \dot{\mathbf{v}} = \mathbf{g}(\mathbf{r}) + \mathbf{a}^d \qquad (1)$$

where $\mathbf{r}$ and $\mathbf{v}$ represent the position and velocity of the rocket, respectively, $\mathbf{g}$ is the gravitational acceleration, and $\mathbf{a}^d$ is the external acceleration. The IIP is the point on the surface of the Earth that meets the free-fall trajectory of the rocket, given the current position and velocity ($\mathbf{r}_0$ and $\mathbf{v}_0$) without any external acceleration ($\mathbf{a}^d = \mathbf{0}$). If the gravitational acceleration is determined based on the Newton's inverse-square law ($\mathbf{g}(\mathbf{r}) = -(\mu / \| \mathbf{r} \|^3) \cdot \mathbf{r}$), the IIP of the rocket in the inertial coordinate frame is computed as

$$\mathbf{i}_p = \frac{\cos(\gamma_0 + \phi)}{\cos \gamma_0} \mathbf{i}_{r0} + \frac{\sin \phi}{\cos \gamma_0} \mathbf{i}_{v0} = [i_{px} \ i_{py} \ i_{pz}]^T \tag{2}$$

where $\mathbf{i}_p$ is the unit vector in the direction of IIP, $\mathbf{i}_{r0}$ and $\mathbf{i}_{v0}$ are the unit vectors for the current position and velocity, respectively, and $\gamma_0$ is the current flight path angle of the rocket. In addition, $\phi$ is the angle of flight (from the current position to the impact point on the surface) obtained by solving the following equation.

$$\frac{r_0}{r_p} = \frac{1 - \cos \phi}{\lambda \cos^2 \gamma_0} + \frac{\cos(\phi + \gamma_0)}{\cos \gamma_0} \tag{3}$$

where $r_p$ is the distance between the center of the Earth and the impact point (= magnitude of the impact point vector), and $\lambda$ is the square of the ratio of $v_0$ and the circular orbital velocity with given $r_0$ ($\equiv (v_0 / v_c)^2 = r_0 v_0^2 / \mu$). The solution of Eq. (3) is written as follows:

$$\sin \phi = \frac{A_1 A_3 + \sqrt{A_1^2 A_3^2 - (A_1^2 + A_2^2)(A_3^2 - A_2^2)}}{(A_1^2 + A_2^2)} \tag{4}$$

where parameters $A_1$, $A_2$, and $A_3$ are defined as

$$A_1 = -\frac{h}{\mu r_0}(\mathbf{r}_0 \cdot \mathbf{v}_0), \ \ A_2 = \frac{h^2}{\mu r_0} - 1, \ \ A_3 = \frac{h^2}{\mu r_p} - 1 \tag{5}$$

In Eq. (5), $h$ ($\equiv r_0 v_0 \cos \gamma_0$) represents the angular momentum of the rocket per unit, which is conserved throughout the orbital flight of the rocket. The IIP latitude and longitude in the inertial coordinate frame can be obtained from elements of $\mathbf{i}_p$ ($i_{px}$, $i_{py}$, and $i_{pz}$) as

$$\text{Lat}_p^I = \arcsin(i_{pz}), \ \ \text{Lon}_p^I = \arctan 2(i_{py}, i_{px}) \tag{6}$$

The IIP should be expressed in rgw Earth-centered Earth-fixed coordinate frame (ECEF) since the impact area also rotates with the Earth. The IIP latitude in the ECEF frame ($\mathrm{Lat}_p^E$) is the same as the value in ECI frame ($\mathrm{Lat}_p^I$), but the compensation in the rotation of the Earth is necessary for IIP longitude ($\mathrm{Lon}_p^E$). The IIP latitude and longitude in the ECEF frame are expressed as

$$\mathrm{Lat}_p^E = \mathrm{Lat}_p^I, \quad \mathrm{Lon}_p^E = \mathrm{Lon}_p^I - \omega_e(t - t_{ref} + t_F) \tag{7}$$

where $\omega_e$ is the Earth's rotational rate, $t$ is the current time, $t_{ref}$ is the reference time, and $t_F$ is the time of flight. The reference time is determined so that the ECI and ECEF frames aligns at $t_{ref}$. The time of flight, $t_F$, is expressed using the result presented by Zarchan (2002) as

$$t_F = \frac{r_0}{v_0 \cos\gamma_0} \left( \frac{\tan\gamma_0(1-\cos\phi)+(1-\lambda)\sin\phi}{(2-\lambda)\left(\dfrac{1-\cos\phi}{\lambda\cos^2\gamma_0}+\dfrac{\cos(\gamma_0+\phi)}{\cos\gamma_0}\right)} + \frac{2\cos\gamma_0}{\lambda(2/\lambda-1)^{1.5}} \arctan\left(\frac{\sqrt{2/\lambda}-1}{\cos\gamma_0\cot(\phi/2)-\sin\gamma_0}\right) \right) \tag{8}$$

### Time Derivatives of Keplerian IIP

The IIP changes when a non-zero external acceleration ($\mathbf{a}^d$) exists as

$$\mathbf{a}^d = a_r^d \mathbf{i}_r + a_\theta^d \mathbf{i}_\theta + a_h^d \mathbf{i}_h \tag{9}$$

where $\mathbf{i}_r (= \mathbf{i}_{r0})$ and $\mathbf{i}_h$ are the unit vectors in the directions of the current position and angular momentum of the rocket, respectively, $\mathbf{i}_\theta$ is defined as $\mathbf{i}_\theta \equiv \mathbf{i}_h \times \mathbf{i}_r$, and $a_r^d / a_h^d / a_\theta^d$ are the components of $\mathbf{a}^d$ in the directions of these three vectors. The time derivatives of the Keplerian IIP introduced in this subsection were presented by Ahn and Roh (2014). In their approach, analytic expres-sions for time derivatives of the IIP were obtained by 1) deriving expressions for $\dot{\phi}$ and $\dot{t}_F$ as functions of $\mathbf{r}_0$, $\mathbf{v}_0$, and $\mathbf{a}^d$; and 2) combining them to determine derivatives of IIP latitude and longitude. The time derivative of flight angle ($\dot{\phi}$) is expressed as

$$\dot{\phi} = -\frac{h}{r_0^2} + D_r^\phi a_r^\phi + D_\theta^\phi a_\theta^\phi,$$

(10)

where $D_r^\phi$ and $D_\theta^\phi$ are defined as

$$D_r^\phi = \frac{h \sin\phi}{\mu(-A_2 \sin\phi + A_1 \cos\phi)}$$

(11)

$$D_\theta^\phi = \frac{2h(r_0 / r_p - \cos\phi) + (\mathbf{r}_0 \cdot \mathbf{v}_0)\sin\phi}{\mu(-A_2 \sin\phi + A_1 \cos\phi)}$$

(12)

Note that the parameters $A_1$, $A_2$, and $A_3$ used in Eqs. (11)-(12) are defined in Eq. (5). Then, the time derivative of the flight tine ($\dot{t}_F$) is expressed as

$$\begin{aligned} \dot{t}_F &= -1 + D_a^t(D_r^a a_r^d + D_\theta^a a_\theta^d) + D_e^t(D_r^e a_r^d + D_\theta^e a_\theta^d) \\ &= -1 + (D_a^t D_r^a + D_e^t D_r^e)a_r^d + (D_a^t D_\theta^a + D_e^t D_\theta^e)a_\theta^d \end{aligned}$$

(13)

where $D_a^t$, $D_e^t$, $D_r^a$, $D_\theta^a$, $D_r^e$, and $D_\theta^e$ are given as

$$\begin{aligned} D_a^t &= \frac{\partial t_F}{\partial a} = \frac{3t_F}{2a} - \frac{1}{a^2 en}\left[\frac{r_p(1 - e\cos E_p)}{\sin E_p} - \frac{r_0(1 - e\cos E_0)}{\sin E_0}\right] \\ D_e^t &= \frac{\partial t_F}{\partial e} = \frac{1}{n}\left[\left(\frac{\cos E_p(1 - e\cos E_p)}{e\sin E_p} - \frac{\cos E_0(1 - e\cos E_0)}{e\sin E_0}\right) - \left(\sin E_p - \sin E_0\right)\right] \end{aligned}$$

(14)

$$D_r^a = \frac{2a^2 v_0 \sin\gamma_0}{\mu}, \;\; D_\theta^a = \frac{2a^2 v_0 \cos\gamma_0}{\mu}$$

(15)

$$D_r^e = \frac{p v_0 \sin\gamma_0}{\mu e}, \; D_\theta^e = \frac{\left(pa - r_0^2\right)v_0 \cos\gamma_0}{\mu ae}$$

(16)

Eqs. (14)-(16) involves orbital parameters of the rocket such as semi-major axis ($a$), eccentricity ($e$), mean motion ($n$), semiparameter ($p$), and current/impact eccentric anomaly values ($E_0/E_p$). All of the parameters can be computed using $\mathbf{r}_0$ and $\mathbf{v}_0$ – the procedure is not presented in this paper (Battin 1999).

To obtain the time derivative of $\mathbf{i}_p$, the angular momentum, $h$, is multiplied to both sides of Eq. (2) to yield:

$$
\begin{aligned}
h \cdot \mathbf{i}_p &= r_0 v_0 \cos(\gamma_0 + \phi) \cdot \mathbf{i}_{r_0} + r_0 v_0 \sin \phi \cdot \mathbf{i}_{v_0} \\
&= r_0 v_0 (\cos \gamma_0 \cos \phi - \sin \gamma_0 \sin \phi) \cdot \mathbf{i}_{r_0} + r_0 \sin \phi \cdot \mathbf{v}_0 \\
&= h \cos \phi \cdot \mathbf{i}_{r_0} - (\mathbf{r}_0 \cdot \mathbf{v}_0) \sin \phi \cdot \mathbf{i}_{r_0} + r_0 \sin \phi \cdot \mathbf{v}_0
\end{aligned}
\tag{17}
$$

By differentiating and rearranging Eq. (17), the expression for $d(\mathbf{i}_p)/dt$ can be obtained as

$$
\frac{d}{dt} \mathbf{i}_p = \begin{pmatrix} \dfrac{di_{px}}{dt} \\[2mm] \dfrac{di_{py}}{dt} \\[2mm] \dfrac{di_{pz}}{dt} \end{pmatrix} = \left( a_r^d \right) \mathbf{d}_r + \left( a_\theta^d \right) \mathbf{d}_\theta + \left( a_h^d \right) \mathbf{d}_h,
\tag{18}
$$

where

$$
\mathbf{d}_r = \frac{1}{h} \Big[ -\big( h \sin \phi + (\mathbf{r}_0 \cdot \mathbf{v}_0) \cos \phi \big) D_r^\phi \cdot \mathbf{i}_{r_0} + \big( r_0 v_0 \cos \phi \big) D_r^\phi \cdot \mathbf{i}_{v_0} \Big],
\tag{19}
$$

$$
\begin{aligned}
\mathbf{d}_\theta &= [ ( r_0 \cos \phi - ( h \sin \phi + (\mathbf{r}_0 \cdot \mathbf{v}_0) \cos \phi ) D_\theta^\phi ) / h ] \cdot \mathbf{i}_{r_0} \\
&\quad + [ ( r_0 v_0 \cos \phi \cdot D_\theta^\phi ) / h ] \cdot \mathbf{i}_{v_0} + [ ( r_0 \sin \phi ) / h ] \cdot \mathbf{i}_\theta + ( -r_0 / h ) \cdot \mathbf{i}_p
\end{aligned}
\tag{20}
$$

$$
\mathbf{d}_h = \frac{1}{h} \big( r_0 \sin \phi \big) \cdot \mathbf{i}_h.
\tag{21}
$$

Time derivatives of IIP latitude and longitude are expressed by differentiating Eqs. (6)-(7) and rearranging them as

$$
\frac{d(\mathrm{Lat}_p^I)}{dt} = \left[ \frac{1}{\cos(\mathrm{Lat}_p^I)} \right] \frac{di_{pz}}{dt}
\tag{22}
$$

$$
\frac{d(\mathrm{Lon}_p^I)}{dt} = \frac{\left[ \cos(\mathrm{Lon}_p^I) \right] \dfrac{di_{py}}{dt} - \left[ \sin(\mathrm{Lon}_p^I) \right] \dfrac{di_{px}}{dt}}{\cos(\mathrm{Lon}_p^I) i_{px} + \sin(\mathrm{Lon}_p^I) i_{py}}.
\tag{23}
$$

$$\frac{d(\text{Lon}_I^E)}{dt} = \frac{d(\text{Lon}_p^I)}{dt} - \omega_e(1 + \dot{t}_F) \qquad (24)$$

***Limitations of the Existing IIP Rate Calculation Procedure***

While the analytic formulation presented in Eqs. (9)-(24) provides the time derivatives of IIP that can be implemented with a relatively tractable procedure, there are some of limitations that should be addressed to improve its applicability. First, the procedure to obtain the inertial IIP rate resorts to mechanical differentiation of IIP definition. Thus, the resultant expression in Eq. (18) is presented without physical interpretation – the meanings of $\mathbf{d}_r$, $\mathbf{d}_\theta$, and $\mathbf{d}_h$ are missing, in particular.

Second, the time derivative vector expressed in a rotating coordinate frame (ECEF coordinate) is not fully developed and presented. While the vector can be reconstructed using the derivatives of IIP latitude and longitude in a rotating frame, it can be expressed by adding a term representing the Earth rate to the inertial IIP rate vector and transforming it with a rotation matrix. Note that typical applications of IIP and its derivatives (e.g. operation of flight safety system (FAA 2000), IIP guidance (Madic 2009)) define the missions on the surface of the Earth – not in an inertial space – and the derivative vectors in a rotating frame are useful.

The new formulation introduced in the next section is developed to address the aforementioned limitations.

## New Formulation for Time Derivative of IIP

Fig. 2 shows the geometry for computing the time derivatives of IIP using the proposed procedure. In the new approach, the IIP change in inertial frame ($\Delta \mathbf{p}$) during a very small time period ($\Delta t$) is separated into two different components: in-plane motion ($\Delta \mathbf{p}_D$) and plane change motion ($\Delta \mathbf{p}_C$). Since the IIP moves on the surface of the Earth, $\Delta \mathbf{p}$ is located on the plane tangential to $\mathbf{p}$ – and is therefore perpendicular to $\mathbf{p}$ ($\Delta \mathbf{p} \perp \mathbf{p}$).

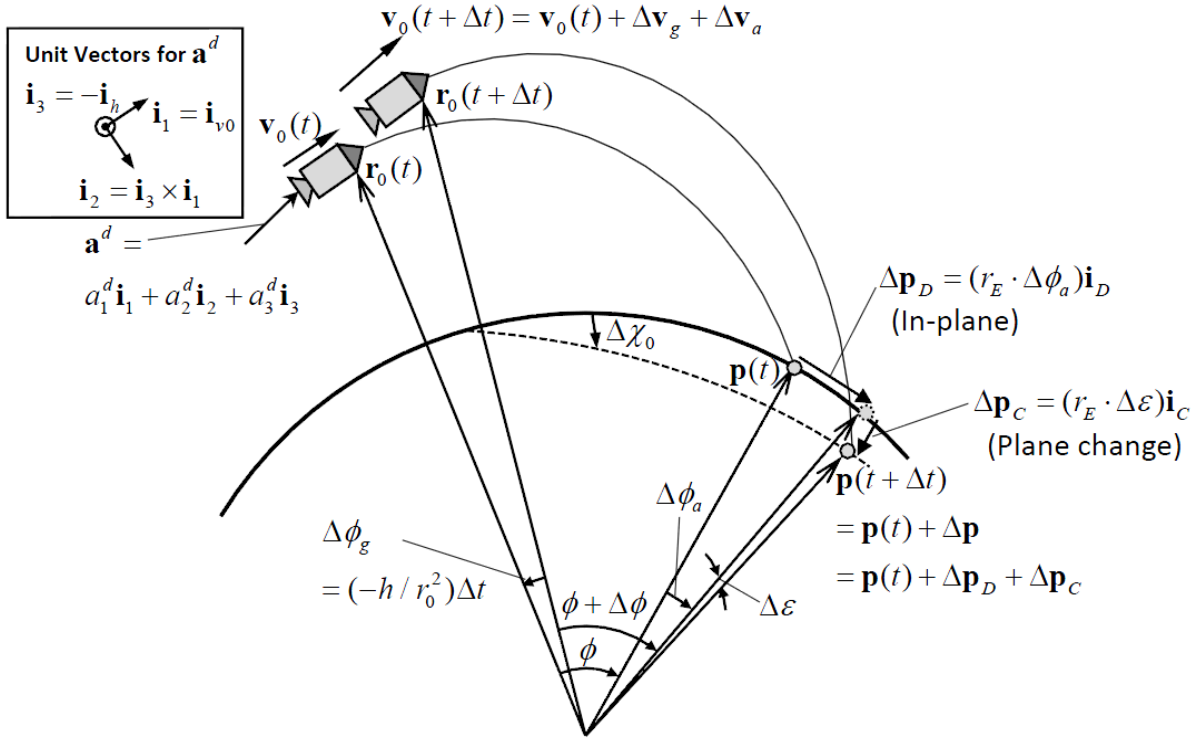

**Fig. 2. Geometry for new procedure to calculate the time derivatives of IIP**

In should be noted that the external acceleration is expressed differently in this new approach as

$$\mathbf{a}^d = a_1^d \mathbf{i}_1 + a_2^d \mathbf{i}_2 + a_3^d \mathbf{i}_3 \qquad (25)$$

where $\mathbf{i}_1$ is parallel to the direction of the velocity, $\mathbf{i}_2$ is the opposite direction of $\mathbf{i}_h$ ($\mathbf{i}_2 = -\mathbf{i}_h$), and $\mathbf{i}_3$ is defined as $\mathbf{i}_3 \equiv \mathbf{i}_1 \times \mathbf{i}_2$. The derivatives of velocity elements due to external acceleration

($\dot{v}_{0,a}$, $\dot{\gamma}_{0,a}$, $\dot{\chi}_{0,a}$) are expressed as:

$$\dot{v}_{0,a} = a_1^d, \ \ \dot{\gamma}_{0,a} = -\frac{a_3^d}{v_0}, \ \ \dot{\chi}_{0,a} = \frac{-a_2^d}{v_0 \cos \gamma_0} \tag{26}$$

The relationship between acceleration elements in Eqs. (9) and (25) is presented in the following equations:

$$a_1^d = \cos \gamma_0 \cdot a_\theta^d + \sin \gamma_0 \cdot a_r^d, \ \ a_2^d = -a_h^d, \ \ a_3^d = \sin \gamma_0 \cdot a_\theta^d - \cos \gamma_0 \cdot a_r^d \tag{27}$$

The following subsections develop formulations for $\Delta \mathbf{p}_D / \Delta \mathbf{p}_C$, and derive the expressions for the IIP derivatives in both the inertial and rotating coordinate frames.

### In-plane Motion

The in-plane motion of the IIP, which occurs in the downrange direction, is discussed first. The unit vector to the downrange direction ($\mathbf{i}_D$) is located in the orbital plane representing the rocket's free fall flight, and can be expressed as a linear combination of $\mathbf{i}_{r0}$ and $\mathbf{i}_{v0}$ as

$$\mathbf{i}_D = \beta_1 \mathbf{i}_{r0} + \beta_2 \mathbf{i}_{v0} \tag{28}$$

The coefficients $\beta_1$ and $\beta_2$ can be found by using the following geometric relationships:

$$\mathbf{i}_D \cdot \mathbf{i}_{r0} = \cos(\pi / 2 - \gamma_0) = \sin \gamma_0 = \beta_1 + (\sin \gamma_0)\beta_2 \tag{29}$$

$$\mathbf{i}_D \cdot \mathbf{i}_{v0} = \cos(\gamma_0 + \phi) = (\sin \gamma_0)\beta_1 + \beta_2 \tag{30}$$

From Eqs. (29)-(30), the coefficients are determined as

$$\beta_1 = -\frac{\sin(\phi + \gamma_0)}{\cos \gamma_0}, \ \ \beta_2 = \frac{\cos \phi}{\cos \gamma_0} \tag{31}$$

The expression for $\Delta \mathbf{p}_D$ is obtained first. The small change in the angle of flight ($\Delta \phi$) is composed of elements contributed by gravity ($\Delta \phi_g$) and external acceleration ($\Delta \phi_a$).

$$\Delta \phi = \Delta \phi_{\mathbf{g}} + \Delta \phi_{\mathbf{a}} \tag{32}$$

Ahn and Roh (2014) pointed out that $\Delta\phi_g\,(=(-h\,/\,r_0^2)\cdot\Delta t$ ) does not result in the change in IIP, and only $\Delta\phi_a$ contributes to $\Delta\mathbf{p}_D$. In particular, the in-plane components of the external acceleration during $\Delta t$ changes $v_0$ and $\gamma_0$ (in-plane velocity components as well), resulting in the change in

$\phi$. Therefore, $\Delta\phi_a$ can be calculated based on the equation.

$$\Delta\phi_{\mathrm{a}}=\frac{\partial\phi}{\partial v_0}\cdot\Delta v_{0,a}+\frac{\partial\phi}{\partial\gamma_0}\cdot\Delta\gamma_{0,a} \tag{33}$$

To obtain the partial derivatives $\partial\phi\,/\,\partial v_0$ and $\partial\phi\,/\,\partial\gamma_0$, Eq. (3) is rearranged using the relationship

$\lambda=r_0 v_0^2\,/\,\mu$ as follows.

$$\left(r_0\,/\,r_p\right)\cos^2\gamma_0=\frac{\mu(1-\cos\phi)}{r_0 v_0^2}+\cos(\phi+\gamma_0)\cdot\cos\gamma_0 \tag{34}$$

The partial derivatives $\partial\phi\,/\,\partial v_0$ and $\partial\phi\,/\,\partial\gamma_0$ can be obtained by differentiating Eq. (34)

$$\frac{\partial\phi}{\partial v_0}=\frac{2\mu(1-\cos\phi)}{r_0 v_0^3}\cdot\frac{1}{D_1} \tag{35}$$

$$\frac{\partial\phi}{\partial\gamma_0}=(\sin(2\gamma_0+\phi)-(r_0\,/\,r_p)\cdot\sin 2\gamma_0)\cdot\frac{1}{D_1} \tag{36}$$

where $D_1$ is defined as

$$D_1\equiv\frac{\mu\sin\phi}{r_0 v_0^2}-\frac{1}{2}[\sin(2\gamma_0+\phi)+\sin\phi] \tag{37}$$

Using Eq. (26), one can express the changes in velocity components ($\Delta v_{0,a}$, $\Delta\gamma_{0,a}$) as

$$\Delta v_{0,a}=a_1^d\Delta t,\ \ \Delta\gamma_0=(-a_3^d\,/\,v_0)\cdot\Delta t \tag{38}$$

Referring to Fig. 2, $\Delta\mathbf{p}_D$ can be expressed as

$$\Delta\mathbf{p}_D=(R_E\Delta\phi_{\mathbf{a}})\mathbf{i}_D=r_p(\frac{\partial\phi}{\partial v_0}\Delta v_{0,a}+\frac{\partial\phi}{\partial\gamma_0}\Delta\gamma_{0,a})\mathbf{i}_D=r_p[(\frac{\partial\phi}{\partial v_0}a_1-\frac{\partial\phi}{\partial\gamma_0}\frac{1}{v_0}a_3)\Delta t]\mathbf{i}_D \tag{39}$$

The following expression for $\dot{\mathbf{p}}_D$ is obtained by taking the limit of Eq. (39) when $\Delta t$ becomes infinitesimally small.

$$\dot{\mathbf{p}}_D = \lim_{\Delta t \to 0} \frac{\Delta \mathbf{p}_D}{\Delta t} = r_p[(\frac{\partial \phi}{\partial v_0}a_1^d - \frac{\partial \phi}{\partial \gamma_0}\frac{1}{v_0}a_3^d)\mathbf{i}_D \qquad (40)$$

*Plane-change Motion*

Plane change motion occurs when the off-plane component of the external acceleration is not zero ($a_2^d \neq 0$). The component changes the azimuth angle ($\Delta\chi_0$), which in turn results in changes in the cross-range angle ($\Delta\varepsilon$) and $\Delta\mathbf{p}_C$. Fig. 3 represents the geometric relationship between the changes in initial azimuth angle ($\Delta\chi_0$) and cross range angle ($\Delta\varepsilon$). Using the law of sines for a spherical triangle **RPP′**, the following relationship can be obtained.

$$\frac{\sin(\pi/2)}{\sin\phi} = \frac{\sin(\Delta\chi_0)}{\sin\Delta\varepsilon_0} \qquad (41)$$

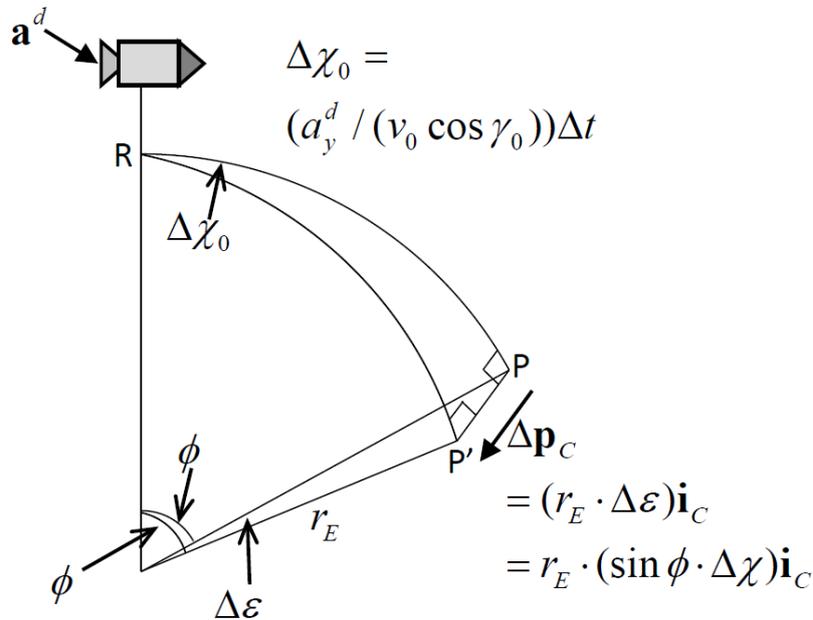

**Fig. 3. Relationship between the change in initial azimuth angle ($\Delta\chi_0$) and the change in cross range angle ($\Delta\varepsilon$)**

Considering that the angles are very small, the approximations for $\sin(\Delta\chi_0)$ and $\sin(\Delta\varepsilon)$ are obtained as $\Delta\chi_0$ and $\Delta\varepsilon$, respectively. Then, the equation relating $\Delta\chi_0$ and $\Delta\varepsilon$ is obtained as

$$\Delta\varepsilon = \sin\phi \cdot \Delta\chi_0 \tag{42}$$

Using Eq. (26), $\Delta\chi_0$ is expressed as

$$\Delta\chi_0 = -\frac{a_2^d}{v_0 \cos\gamma_0} \cdot \Delta t \tag{43}$$

Then $\Delta\mathbf{p}_C$ is expressed as the following equation.

$$\Delta\mathbf{p}_C = (R_E\Delta\varepsilon)\mathbf{i}_C = (-\frac{\sin\phi}{v_0\cos\gamma_0}a_2^d)(\Delta t)\mathbf{i}_C \tag{44}$$

The plane change IIP derivative, $\dot{\mathbf{p}}_C$, is obtained by taking the limit of $\Delta\mathbf{p}_C / \Delta t$ as

$$\dot{\mathbf{p}}_C = \lim_{\Delta t \to 0}\frac{\Delta\mathbf{p}_C}{\Delta t} = (-\frac{\sin\phi}{v_0\cos\gamma_0}a_2^d)\mathbf{i}_C \tag{45}$$

### *IIP Derivatives in ECI/ECEF Coordinate Frames*

Fig. 4 shows the geometry in a plane tangential to $\mathbf{p}$ representing the components of IIP derivative with respect to inertial and rotational frames. The formulations in the previous subsections are expressed in the $X_I Y_I Z_I$ coordinate system. The overall derivative of inertial IIP ($\dot{\mathbf{p}}$) is the sum of the $\dot{\mathbf{p}}_D$ and $\dot{\mathbf{p}}_C$.

$$\dot{\mathbf{p}} = r_p\mathbf{i}_p = \dot{\mathbf{p}}_D + \dot{\mathbf{p}}_C \tag{46}$$

The derivatives of the IIP latitude/longitudes can be obtained using the following equations:

$$\frac{d(\mathrm{Lat}_p^I)}{dt} = (\dot{\mathbf{p}} \cdot \mathbf{i}_N) / r_p \tag{47}$$

$$\frac{d(\mathrm{Lon}_p^I)}{dt} = (\dot{\mathbf{p}} \cdot \mathbf{i}_E) / [r_p \cos(\mathrm{Lat}_p^I)] \tag{48}$$

where $\mathbf{i}_E$ and $\mathbf{i}_N$ are defined using components of $\mathbf{i}_p$ (defined in Eq. (2)) as:

$$\mathbf{i}_E = \frac{1}{\sqrt{i_{px}^2 + i_{py}^2}} \cdot \begin{bmatrix} -i_{py} \\ i_{px} \\ 0 \end{bmatrix} = -\mathbf{i}_W \qquad (49)$$

$$\mathbf{i}_N = \mathbf{i}_p \times \mathbf{i}_E \qquad (50)$$

Time derivatives of $\mathrm{Lon}_p^E$ ( $d(\mathrm{Lon}_p^E)/dt$ ) can be computed using Eq. (24).

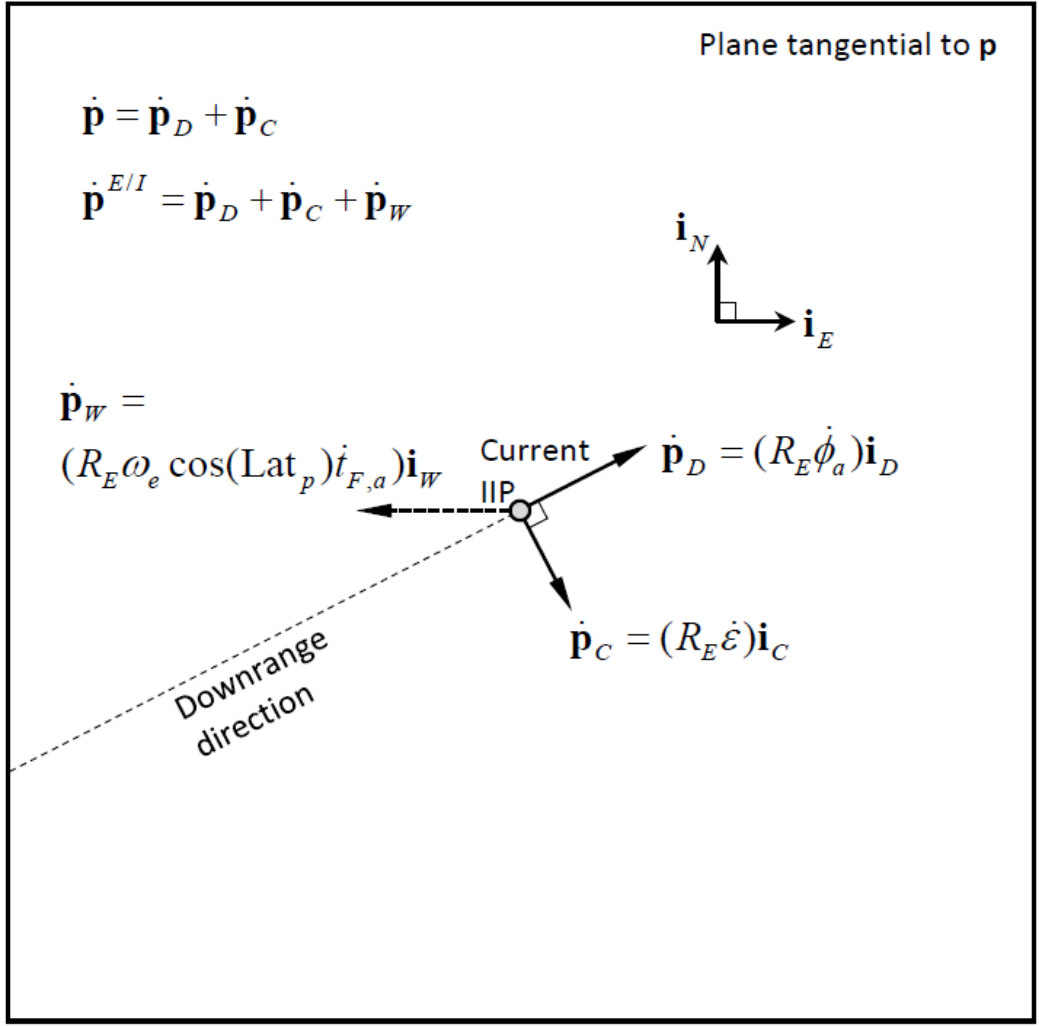

**Fig. 4. Calculation of IIP derivatives in inertial and rotational frames (expressed in $X_I Y_I Z_I$ coordinate system)**

Finally, the formulation for the IIP derivative in ECEF coordinate frame is presented. To compensate for the effects of the Earth rotation, $\dot{\mathbf{p}}^{E/I}$ is expressed as:

$$\dot{\mathbf{p}}^{E/I} = \dot{\mathbf{p}} + \dot{\mathbf{p}}_W \qquad (51)$$

where $\dot{\mathbf{p}}_W$ is defined as

$$\dot{\mathbf{p}}_W = (R_E \omega_E \cos(\text{Lat}_p^I))(1 + \dot{t}_F)\mathbf{i}_W = (R_E \omega_E \cos(\text{Lat}_p^I))(\dot{t}_{F,a})\mathbf{i}_W \tag{52}$$

In Eq. (52), $\dot{t}_{F,a}(=1+\dot{t}_F)$ is the derivative of the time of flight attributed to the external acceleration and $\mathbf{i}_W$ is defined in Eq. (49). Note that $\dot{\mathbf{p}}^{E/I}$ is the IIP derivative in a rotating frame expressed in $X_I Y_I Z_I$ coordinates. To express the IIP derivative in a rotating frame in $X_E Y_E Z_E$ coordinates, $\dot{\mathbf{p}}^{E/I}$ should be transformed using the transformation matrix $\mathrm{T}_I^E$ as

$$\dot{\mathbf{p}}^E = \mathrm{T}_I^E \dot{\mathbf{p}}^{E/I} \tag{53}$$

where $\mathrm{T}_I^E$ is defined as

$$\mathrm{T}_I^E = \begin{bmatrix} \cos \omega_e t_F & \sin \omega_e t_F & 0 \\ -\sin \omega_e t_F & \cos \omega_e t_F & 0 \\ 0 & 0 & 1 \end{bmatrix} \tag{54}$$

### *Verification of the Proposed Formulation through Numerical Simulation*

Numerical simulations were conducted to verify the proposed formulation by comparison with existing IIP derivative expressions. The trajectory of a two-stage launch vehicle introduced in Ahn and Roh (2014) was generated and used for comparing the time derivatives of IIP obtained using the two different formulations. Table 1 summarizes the launch vehicle information used for the numerical simulation, and Fig. 5 presents the profiles of external acceleration components ($a_r^d$, $a_\theta^d$, and $a_h^d$) applied to the vehicle during its flight.

**Table 1. Launch vehicle configuration for the verification simulation** (Ahn and Roh (2014))

| Parameter | First Stage | Payload Fairing | Second Stage | Satellite |
|---|---|---|---|---|
| Structural mass, kg | 7,000 | 300 | 450 | 250 |
| Propellant mass, kg | 100,000 | | 550 | |
| Thrust, $10^3$ kgf | 32-150 | | 2.8 | |
| Specific impulse, s | 320 | | 280 | |
| Burning time, s | 343 | | 55 | |

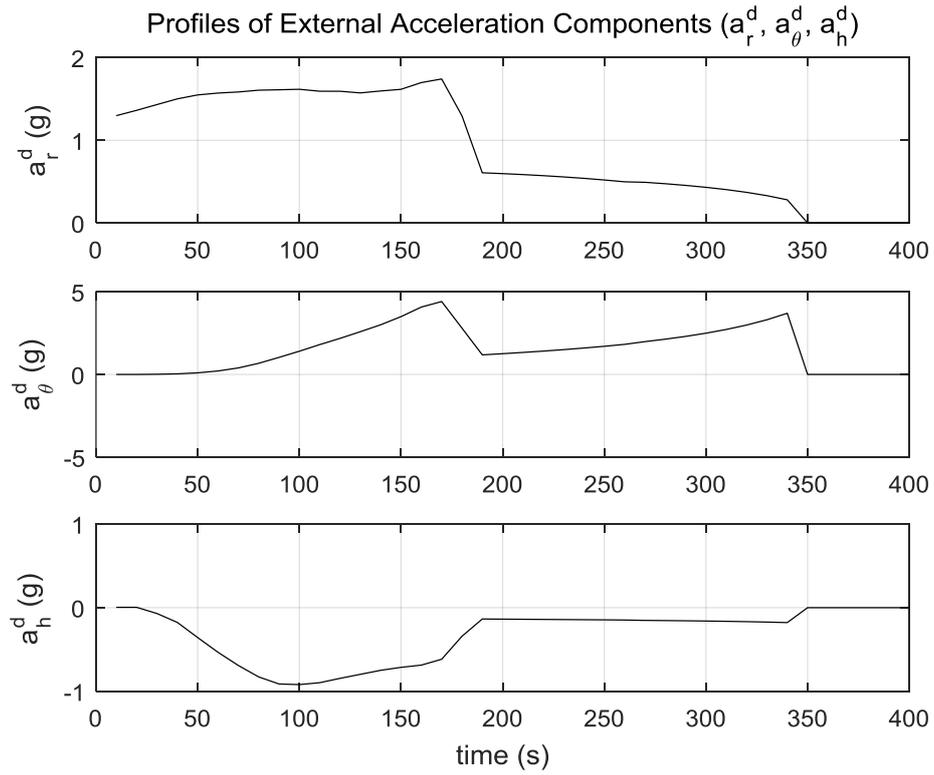

**Fig. 5. External acceleration profile for numerical simulation**

Figs. 6 and 7 compare the IIP latitude and longitude profiles ($d(\text{Lat}_p^I)/dt$ and $d(\text{Lon}_p^E)/dt$) and the components of the IIP unit vector ($i_{px}, i_{py}$, and $i_{pz}$) computed using two different formulations. It can be easily observed that the time derivative values obtained by different formulations were exactly matched. This indicates that the procedure proposed in this paper (presented as Eqs. (25)-(54)) yields the same results as the existing algorithm (presented as Eqs. (9)-(24)) developed by Ahn and Roh (2014).

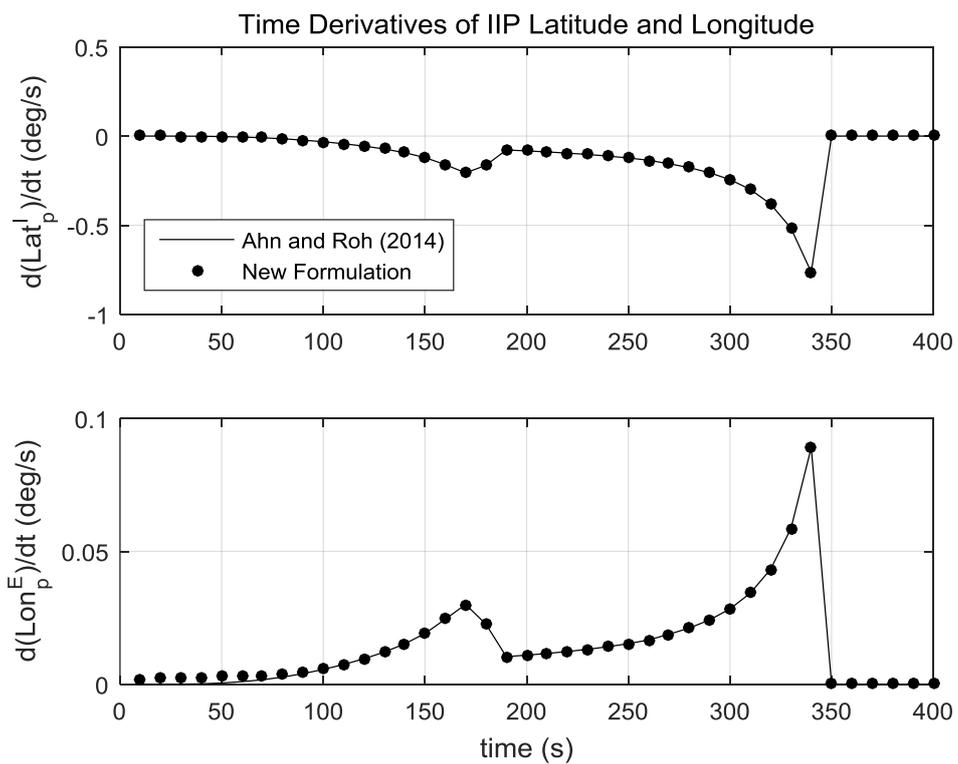

**Fig. 6. Comparison of IIP latitude and longitude profiles**

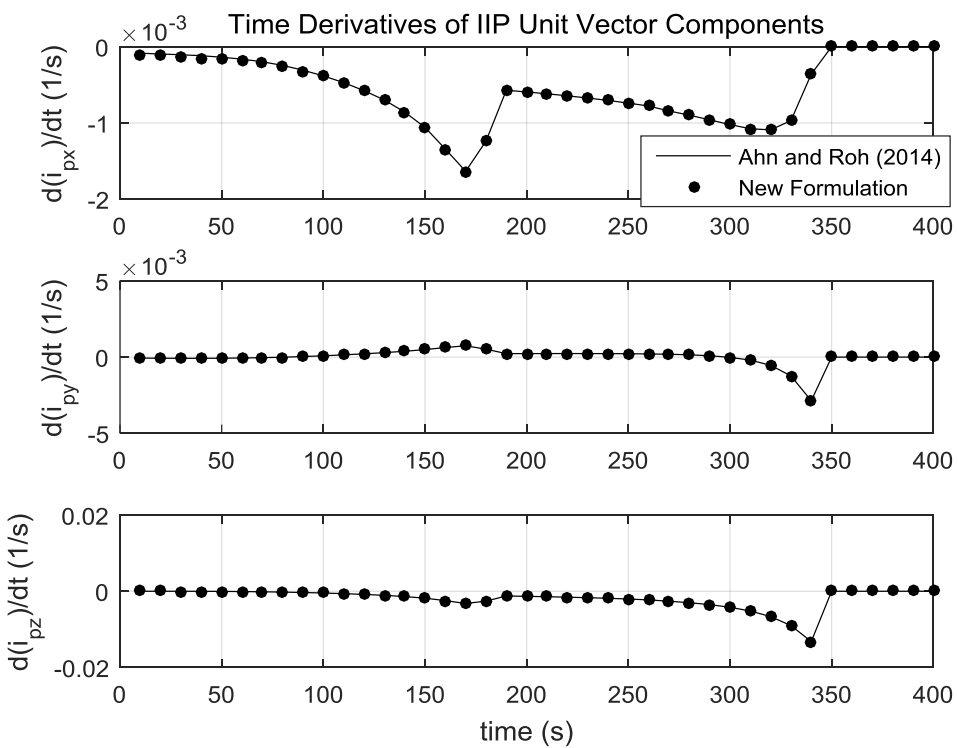

**Fig. 7. Comparison of IIP unit vector component profiles**

*Discussion: Advantages of the New Formulation*

The proposed procedure successfully addresses the drawbacks of the existing formulation mentioned in the previous section. The expression for the inertial IIP derivative presented in Eq. (46) is decomposed into two components having physical interpretations: the downrange component ($\dot{\mathbf{p}}_D$) representing the in-plane motion, and the crossrange component ($\dot{\mathbf{p}}_C$) representing the plane-change motion. This clear decomposition was missing in Eq. (18) – the equivalent expression in the existing formulation. Note that the dependency of basis vectors $\mathbf{d}_r$ and $\mathbf{d}_\theta$ in Eq. (18) were observed through numerical simulation ($\mathbf{d}_r \times \mathbf{d}_\theta = \mathbf{0}$), while its mathematical proof is not presented in this paper. The proposed method also enables simpler implementation than the existing formulation by providing a streamlined procedure – particularly the time derivatives of IIP in a rotational frame ($\dot{\mathbf{p}}^E$), which requires complex vector reconstruction using latitude/longitude IIP rates. In summary, the new formulation can be used as a replacement for the existing formulation without degradation of results and with improved computational efficiency and understandability.

# Conclusions

A new analytic formulation to express the time derivatives of the IIP of a rocket is proposed that utilizes the IIP kinematics on the surface of the Earth. In the proposed procedure, the downrange direction (in-plane motion) and crossrange direction (plane change motion) derivatives are computed separately and then combined to express the overall IIP rate in a systematic and physically explainable way. The resultant expressions obtained by the proposed procedure were verified by comparison with the existing formulation along the trajectory of a launch vehicle obtained by numerical simulation. The main advantages of the proposed procedure are the improved understandability and ease of implementation, both of which justify its replacement for the existing formulation.

## Acknowledgements

This work was prepared at the Korea Advanced Institute of Science and Technology, Department of Aerospace Engineering, under a research grant from the National Research Foundation of Korea (NRF-2013M1A3A3A02042461). Authors thank the National Research Foundation of Korea for the support of this work.

## Notation

*The following symbols are used in this paper*:

| | |
|---|---|
| $\mathbf{r}_0, \mathbf{i}_{r_0}, r_0$ | = current position vector, its unit vector, and its magnitude |
| $\mathbf{v}_0, \mathbf{i}_{v_0}, v_0$ | = current velocity vector, its unit vector, and its magnitude |
| $\mathbf{r}_p, \mathbf{i}_p, r_p$ | = instantaneous impact point vector, its unit vector, and its magnitude |
| $\mathbf{h}, \mathbf{i}_h, h$ | = angular momentum vector, its unit vector, and its magnitude |
| $\mathbf{i}_\theta$ | = unit vector of tangential direction |
| $\mathbf{i}_D, \mathbf{i}_C$ | = unit vectors of downrange and crossrange directions |
| $\mathbf{i}_N, \mathbf{i}_E$ | = unit vector of North and East directions |
| $\gamma_0, \ \chi_0$ | = current inertial flight path and azimuth angles |
| $\mu$ | = gravitational constant of the Earth |
| $\phi$ | = angle of flight of a rocket |
| $\varepsilon$ | = cross range angle of a rocket |
| $t_F$ | = flight time of a rocket |
| $\mathbf{a}^d, a_r^d, a_\theta^d, a_h^d$ | = external acceleration vector and its components in radial, tangential, and linear momentum directions |
| $a_1^d, a_2^d, a_3^d$ | = components of external acceleration vector in proposed formulation |
| $v_c$ | = circular orbital speed |
| $\lambda$ | = squared ratio of current speed to the circular orbital speed |
| $a$ | = semi-major axis of an orbit |

$e$ $\qquad$ = eccentricity of an orbit

$n$ $\qquad$ = mean motion of an orbit

$E_0, E_p$ $\qquad$ = current and impact eccentric anomalies of a rocket

$M_0, M_p$ $\qquad$ = current and impact mean anomalies of a rocket

$\text{Lat}_p^I$ $\qquad$ = impact point latitude

$\text{Lon}_p^I, \text{Lon}_p^E$ $\qquad$ = impact point longitudes in inertial and Earth-fixed frame

$D_r^\phi, D_\theta^\phi$ $\qquad$ = flight angle derivative sensitivities with respect to the radial and tangential disturbing acceleration components

$D_r^t, D_\theta^t$ $\qquad$ = flight time derivative sensitivities with respect to the radial and tangential disturbing acceleration components

$\mathbf{d}_r, \mathbf{d}_\theta, \mathbf{d}_h$ $\qquad$ = IIP vector derivative sensitivities with respect to the radial, tangential, and off-plane disturbing acceleration components

$\dot{\mathbf{p}}, \dot{\mathbf{p}}_D, \dot{\mathbf{p}}_C, \dot{\mathbf{p}}_W$ $\qquad$ = total time derivative of IIP and its decomposed elements caused by changes in downrange, crossrange, and time of flight